\documentclass[12pt]{article}
\usepackage{amsmath,amsthm,amsfonts,amssymb,amscd}
\newtheorem {Lemma}{Lemma}[section]
\newtheorem {Theorem}{Theorem}[section]

\newenvironment {Proof} {\noindent {\bf Proof.}}{\quad $\square$\par\vspace{3mm}}

\allowdisplaybreaks [4]

= 0.0in \topmargin = 0.3in \oddsidemargin=0.1in
\begin{document}

\title{Bounds for the spectral radius of nonnegative matrices}

\author{Rundan Xing, Bo Zhou\footnote{Corresponding author.}
\\
Department of Mathematics, South China Normal University,\\
Guangzhou 510631, P. R. China\\
email: {\tt zhoubo@scnu.edu.cn}}

\date{}
\maketitle

\begin{abstract}
We give upper and lower bounds for the spectral radius of a
nonnegative matrix by using its average $2$-row sums, and
characterize the equality cases if the matrix is irreducible. We
also apply these bounds to various nonnegative matrices associated
with a graph, including the adjacency matrix, the signless Laplacian
matrix, the distance matrix, the distance signless Laplacian matrix,
and the reciprocal distance matrix. 
\end{abstract}

\section{Introduction}

Let $A=(a_{ij})$ be an $n\times n$ nonnegative matrix. The spectral
radius of $A$, denoted by $\rho(A)$, is the largest modulus of
eigenvalues of $A$. Moreover, if $A$ is symmetric, then $\rho(A)$ is
equal to the largest eigenvalue of $A$. See~\cite{BP,HJ,Mi} for some
well-known properties of the spectral radius of nonnegative
matrices.

We consider simple graphs.  Let $G$ be a graph with vertex set
$V(G)=\{v_1,\dots,v_n\}$ and edge set $E(G)$.

The adjacency matrix of $G$ is the $n\times n$ matrix
$A(G)=(a_{ij})$, where $a_{ij}=1$ if $v_i$ and $v_j$ are adjacent in
$G$, and $0$ otherwise~\cite{CDS}. The spectral radius of the
adjacency matrix of graphs has been studied extensively, see
\cite{Ho,LW,SW}.

For $1\le i\le n$, let $d_i$ be the degree of vertex $v_i$, which is
the number of vertices adjacent to $v_i$, in $G$. Let $\Delta(G)$ be
the degree diagonal matrix $\text{diag}(d_1,\dots,d_n)$. The
signless Laplacian matrix of $G$ is the $n\times n$ matrix
$Q(G)=(q_{ij})=\Delta(G)+A(G)$~\cite{CRS}. The spectral radius of
the signless Laplacian matrix of graphs has received much attention
recently, see~\cite{Yu}.

Suppose that $G$ is connected.

The distance matrix of $G$ is the $n\times n$ matrix
$D(G)=(d_{ij})$, where $d_{ij}$ is  the distance between vertices
$v_i$ and $v_j$, which is the length of a shortest path connecting
them, in $G$~\cite{JT}.

For $1\le i\le n$, let $D_i$ be the transmission of vertex $v_i$ in
$G$, which is the sum of distances between $v_i$ and (other)
vertices of $G$. Obviously, $D_i=r_i(D(G))$ for $1\le i\le n$. Let
$Tr(G)$ be the transmission diagonal matrix
$\text{diag}(D_1,\dots,D_n)$. The distance signless Laplacian matrix
of $G$ is the $n\times n$ matrix
$DQ(G)=(dq_{ij})=Tr(G)+D(G)$~\cite{AH}.

The reciprocal distance matrix (also called the Harary matrix) of
$G$ is the $n\times n$ matrix $R(G)=(r_{ij})$, where
$r_{ij}=\frac{1}{d_{ij}}$ if $i\ne j$, and $r_{ii}=0$  for $1\le
i\le n$~\cite{JT}.

There are also some results on the  spectral radius of the distance
matrix and some distance-based matrices of connected graphs,
see~\cite{Zhou}.

For $1\le i\le n$, the $i$-th row sum of the nonnegative matrix
$A=(a_{ij})_{n\times n}$ is $r_i(A)=\sum_{j=1}^na_{ij}$. Very
recently, Duan and Zhou~\cite{DZ} found upper and lower bounds for
the spectral radius of a nonnegative matrix using its row sums, and
characterized the extremal cases if the matrix is irreducible. They
also applied the bounds to the matrices associated with a graph as
mentioned above.

In the whole of this paper, suppose that $r_i(A)>0$ for each $1\le
i\le n$. The $i$-th average $2$-row sum of $A$ is defined as
\[
m_i(A)=\frac{\sum_{k=1}^{n}a_{ik}r_k(A)}{r_i(A)}=a_{ii}+\frac{\sum_{1\le
k\le n\atop k\ne i}a_{ik}r_k(A)}{r_i(A)}.
\]
For an $n$-vertex graph $G$ which contains no isolated vertices,
$m_i(A(G))=\frac{\sum_{v_iv_j\in E(G)}d_j}{d_i}$ with $1\le i\le n$,
which is known as the average $2$-degree of vertex $v_i$ in
$G$~\cite{Ca,Me}. Chen et al.~\cite{CPZ} gave upper bound for the
spectral radius of the adjacency matrix of a connected graph using
the two largest average $2$-degrees, which was refined very recently
by Huang and Weng~\cite{HW} using average $2$-degrees.

Motivated by the work of~\cite{DZ,HW}, we give upper and lower
bounds for the spectral radius of a nonnegative matrix by using its
average $2$-row sums, and characterize the equality cases if the
matrix is irreducible. We also apply these results to various
matrices associated with a graph as mentioned above.

\section{Bounds for the spectral radius of a nonnegative matrix}

The following lemma is well known.

\begin{Lemma} \label{r-min-max}\cite{Mi}
If $A$ is an $n\times n$  nonnegative   matrix, then
\[
\min_{1\le i\le n}r_i(A)\le\rho(A)\le\max_{1\le i\le n}r_i(A).
\]
Moreover, if $A$ is irreducible, then either equality holds if and
only if $r_1(A)=\dots=r_n(A)$.
\end{Lemma}

The following lemma is the starting point of this paper, which has
been given in~\cite{ZL} for an irreducible nonnegative matrix. We
include a proof here for completeness.

\begin{Lemma}\label{m1}
Let $A$ be an $n\times n$ nonnegative matrix. Then
\[
\min_{1\le i\le n}m_i(A)\le\rho(A)\le\max_{1\le i\le n}m_i(A).
\]
Moreover, if $A$ is irreducible, then either equality holds if and
only if $m_1(A)=\dots=m_n(A)$.
\end{Lemma}

\begin{Proof}
Let $A=(a_{ij})$. For $U=\text{diag}(r_1(A),\dots,r_n(A))$, let
$B=(b_{ij})=U^{-1}AU$. Evidently,
$b_{ij}=\frac{a_{ij}r_j(A)}{r_i(A)}$ for $1\le i,j\le n$. Then
$r_i(B)=\sum_{k=1}^n b_{ik}=\frac{\sum_{k=1}^n
a_{ik}r_k(A)}{r_i(A)}=m_i(A)$ for $1\le i\le n$. The result follows
easily from Lemma~\ref{r-min-max}.
\end{Proof}

\begin{Theorem} \label{nonnega-up}
Let $A=(a_{ij})$ be an $n\times n$ nonnegative matrix with average
$2$-row sums $m_1,\dots,m_n$, where $m_1\ge\dots\ge m_n$. Let $M$ be
the largest diagonal element, and $N$ the largest off-diagonal
element of $A$. Suppose that $N>0$. Let
$b=\max\left\{\frac{r_j(A)}{r_i(A)}:1\le i,j\le n\right\}$. For
$1\le l\le n$, let
\[
\phi_l=\frac{m_{l}+M-Nb+
\sqrt{(m_{l}-M+Nb)^2+4Nb\sum_{i=1}^{l-1}(m_i-m_l)}}{2}.
\]
Then $\rho(A)\le\phi_l$ for $1\le l\le n$. Moreover, if $A$ is
irreducible, then $\rho(A)=\phi_l$ if and only if $m_1=\dots=m_n$ or
for some $2\le t\le l$, $A$ satisfies the following conditions:

(i) $a_{ii}=M$ for $1\le i\le t-1$,

(ii) $m_t=\dots=m_n$,

(iii) $a_{ik}=N$ and $\frac{r_k(A)}{r_i(A)}=b$ for $1\le i\le n$,
$1\le k\le t-1$ and $k\ne i$.
\end{Theorem}

\begin{Proof}
For convenience, let $r_i=r_i(A)$ for $1\le i\le n$. Suppose that
$a_{ii}=M$ for some $1\le i\le m$. Then
\[
m_1\ge m_i=a_{ii}+\frac{\sum_{1\le k\le n\atop k\ne
i}a_{ik}r_k}{r_i}\ge a_{ii}=M.
\]
If $l=1$, then $\phi_l=\frac{m_1+M-Nb+|m_1-M+Nb|}{2}=m_1$, and thus
the result follows immediately from Lemma~\ref{m1}. Suppose in the
following that $2\le l\le n$.

Let $U=\text{diag}(r_1x_1,\dots,r_{l-1}x_{l-1},r_l,\dots,r_n)$,
where $x_i\ge 1$ is a variable to be determined later for $1\le i\le
l-1$. Let $B=U^{-1}AU$. Obviously, $A$ and $B$ are unitary similar,
and thus have the same eigenvalues. Recall that $m_i=\sum_{k=1}^n
a_{ik}\frac{r_k}{r_i}$ for $1\le i\le n$. For $1\le i\le l-1$, since
$a_{ii}\le M$, and $a_{ik}\le N$, $\frac{r_k}{r_i}\le b$ for $k\ne
i$, we have
\begin{eqnarray*}
   r_i(B)
&=&
   \frac{1}{x_i}\left(\sum_{k=1}^{l-1}a_{ik}\frac{r_k}{r_i}x_k+\sum_{k=l}^n a_{ik}\frac{r_k}{r_i}\right)\\
&=&
   \frac{1}{x_i}\left(\sum_{k=1}^{l-1}a_{ik}\frac{r_k}{r_i}(x_k-1)+\sum_{k=1}^n
   a_{ik}\frac{r_k}{r_i}\right)\\
&=&
   \frac{1}{x_i}\left(\sum_{1\le k\le l-1\atop k\ne
   i}a_{ik}\frac{r_k}{r_i}(x_k-1)+a_{ii}(x_i-1)+m_i\right)\\
&\le&
   \frac{1}{x_i}\left(Nb\sum_{1\le k\le l-1\atop k\ne
   i}(x_k-1)+M(x_i-1)+m_i\right)\\
&=&
   \frac{1}{x_i}\left(Nb\sum_{k=1}^{l-1}(x_k-1)+(M-Nb)(x_i-1)+m_i\right)
\end{eqnarray*}
with equality if and only if (a) and (b) hold: (a) $x_i=1$ or
$a_{ii}=M$, (b) $x_k=1$ or  $a_{ik}=N$, $\frac{r_k}{r_i}=b$, where
$1\le k\le l-1$ and $k\ne i$. For $l\le i\le n$, since $m_i\le m_l$,
and $a_{ik}\le N$, $\frac{r_k}{r_i}\le b$ for $1\le k\le l-1$, we
have
\begin{eqnarray*}
   r_i(B)
&=&
   \sum_{k=1}^{l-1}a_{ik}\frac{r_k}{r_i}x_{k}+\sum_{k=l}^n a_{ik}\frac{r_k}{r_i}\\
&=&
   \sum_{k=1}^{l-1}a_{ik}\frac{r_k}{r_i}(x_{k}-1)+\sum_{k=1}^n a_{ik}\frac{r_k}{r_i}\\
&=&
   \sum_{k=1}^{l-1}a_{ik}\frac{r_k}{r_i}(x_{k}-1)+m_i\\
&\le&
   Nb\sum_{k=1}^{l-1}(x_{k}-1)+m_l
\end{eqnarray*}
with equality if and only if (c) and (d) hold: (c) $m_i=m_l$, (d)
$x_k=1$ or $a_{ik}=N$, $\frac{r_k}{r_i}=b$, where $1\le k\le l-1$.

Recall that for $1\le l\le n$,
\[
\phi_l=\frac{m_{l}+M-Nb+
\sqrt{(m_{l}-M+Nb)^2+4Nb\sum_{k=1}^{l-1}(m_k-m_l)}}{2},
\]
and thus
$\phi_l^2-(m_l+M-Nb)\phi_l+m_l(M-Nb)-Nb\sum_{k=1}^{l-1}(m_k-m_l)=0$,
i.e.,
\[
Nb\sum_{k=1}^{l-1}(m_k-m_l)=(\phi_l-m_l)(\phi_l-M+Nb).
\]
Note that $Nb>0$. If $\sum_{k=1}^{l-1}(m_k-m_l)>0$, then
$\phi_l>\frac{m_l+M-Nb+|m_l-M+Nb|}{2}\ge\frac{m_l+M-Nb-(m_l-M+Nb)}{2}=M-Nb$,
and if $m_1=\dots=m_l$, then since $m_1\ge M$, we have
$\phi_l=\frac{m_1+M-Nb+|m_1-M+Nb|}{2}>\frac{m_1+M-Nb-(m_1-M+Nb)}{2}=M-Nb$.
It follows that $\phi_l-M+Nb>0$. For $1\le i\le l-1$, let
$x_i=1+\frac{m_i-m_l}{\phi_l-M+Nb}$. Obviously, $x_i\ge 1$, and
\[
Nb\sum_{k=1}^{l-1}(x_k-1)=\frac{Nb\sum_{k=1}^{l-1}(m_k-m_l)}{\phi_l-M+Nb}
=\phi_l-m_l.
\]
Thus for $1\le i\le l-1$,
\begin{eqnarray*}
   r_i(B)
&\le&
   \frac{1}{x_i}\left(Nb\sum_{k=1}^{l-1}(x_k-1)+(M-Nb)(x_i-1)+m_i\right)\\
&=&
   \frac{(\phi_l-m_l)+(M-Nb)\cdot\frac{m_i-m_l}{\phi_l-M+Nb}+m_i}{1+\frac{m_i-m_l}{\phi_l-M+Nb}}\\
&=&
   \phi_l,
\end{eqnarray*}
and for $l\le i\le n$,
\[
r_i(B)\le Nb\sum_{k=1}^{l-1}(x_{k}-1)+m_l=(\phi_l-m_l)+m_l=\phi_l.
\]
Hence by Lemma~\ref{r-min-max}, $\rho(A)= \rho(B)\le \max_{1\le i\le
n}r_i(B)\le\phi_l$.

Now suppose that $A$ is irreducible. Then $B$ is also irreducible.

Suppose that $\rho(A)=\phi_l$ for some $2\le l\le n$. Then
$\rho(B)=\max_{1\le i\le n}r_i(B)=\phi_l$, which, by
Lemma~\ref{r-min-max}, implies that $r_1(B)=\dots=r_n(B)=\phi_l$,
and thus from the above arguments, (a) and (b) hold for $1\le i\le
l-1$, and (c) and (d) hold for $l\le i\le n$. If $m_1=m_l$, then
since from (c), $m_i=m_l$ for $l\le i\le n$, we have
$m_1=\dots=m_n$. Now assume that $m_1>m_l$. Let $t$ be the smallest
integer such that $m_t=m_l$, where $2\le t\le l$. From (c), we now
have $m_t=\dots=m_l=\dots=m_n$, implying that (ii) holds. For $1\le
i\le t-1$, since $m_i>m_l$, we have $x_i>1$. Now (i) and (iii)
follow from (a), (b) for $1\le i\le l-1$ and (d) for $l\le i\le n$.

Conversely, if $m_1=\dots=m_n$, then
$\phi_l=\frac{m_1+M-Nb+|m_1-M+Nb|}{2}=m_1$, and thus by
Lemma~\ref{m1}, $\rho(A)=m_1=\phi_l$. If (i)--(iii) hold, then (a)
and (b) hold for $1\le i\le l-1$, and (c) and (d) hold for $l\le
i\le n$, implying that $r_i(B)=\phi_l$ for $1\le i\le n$, and thus
by Lemma~\ref{r-min-max}, $\rho(A)=\rho(B)=\phi_l$.
\end{Proof}

Let $I_n$ and $J_n$ be the $n\times n$ identity matrix and the
$n\times n$ all-one matrix, respectively.

Under the conditions of Theorem~\ref{nonnega-up}, we have
$\phi_l\ge\phi_{l+1}$ if and only if
\begin{eqnarray*}
& &
   m_l-m_{l+1}+\sqrt{(m_l-M+Nb)^2+4Nb\sum_{i=1}^{l-1}(m_i-m_l)}\\
&\ge&
   \sqrt{(m_{l+1}-M+Nb)^2+4Nb\sum_{i=1}^l(m_i-m_{l+1})},
\end{eqnarray*}
which is equivalent to
\begin{eqnarray*}
& &
   (m_l-m_{l+1})\sqrt{(m_l-M+Nb)^2+4Nb\sum_{i=1}^{l-1}(m_i-m_l)}\\
&\ge&
   (m_l-m_{l+1})(2Nbl+M-Nb-m_l).
\end{eqnarray*}
Note that
\[\sqrt{(m_l-M+Nb)^2+4Nb\sum_{i=1}^{l-1}(m_i-m_l)}\ge 2Nbl+M-Nb-m_l
\]
if and only if $\sum_{i=1}^lm_i\ge l(Nbl+M-Nb)$. Thus if
$\sum_{i=1}^l m_i\ge l(Nbl+M-Nb)$, then $\phi_l\ge\phi_{l+1}$, and
if $\sum_{i=1}^l m_i<l(Nbl+M-Nb)$, then $\phi_l\le\phi_{l+1}$.

For $1\le i\le n$, since $a_{ii}\le M$, and $a_{ij}\le N$,
$\frac{r_j(A)}{r_i(A)}\le b$ for $j\ne i$, we have
\[
\sum_{i=1}^n m_i=\sum_{i=1}^n\left(a_{ii}+\sum_{1\le j\le n\atop
j\ne i}a_{ij}\frac{r_j(A)}{r_i(A)}\right)\le\sum_{i=1}^n(M+Nb(n-1))=
n(Nbn+M-Nb)
\]
with equality if and only if $a_{ii}=M$ for $1\le i\le n$, and
$a_{ij}=N$, $\frac{r_j(A)}{r_i(A)}=b$ for $1\le i,j\le n$ and $i\ne
j$, i.e., $A=MI_n+(N-M)J_n$. Suppose that $A\ne MI_n+(N-M)J_n$,
which implies that $\sum_{i=1}^n m_i<n(Nbn+M-Nb)$. Note that $m_1\ge
M$. Let $l$ be the smallest integer such that $\sum_{i=1}^l
m_i<l(Nbl+M-Nb)$, where $2\le l\le n$. For $1\le s\le l-1$, we have
by the choice of $l$ that $\sum_{i=1}^sm_i\ge s(Nbs+M-Nb)$, and thus
$\phi_1\ge\dots\ge\phi_{l}$. For $l\le s\le n$, we are to show that
$\sum_{i=1}^sm_i< s(Nbs+M-Nb)$ by induction on $s$. The case $s=l$
has been done from our choice of $l$. Suppose that $\sum_{i=1}^s
m_i<s(Nbs+M-Nb)$ for some $l\le s\le n-1$. Then $m_s<Nbs+M-Nb$,
which, together with the fact that $m_{s+1}\le m_s$, implies that
$\sum_{i=1}^{s+1}m_i< s(Nbs+M-Nb)+(Nbs+M-Nb)<(s+1)(Nb(s+1)+M-Nb)$.
It follows that $\sum_{i=1}^sm_i< s(Nbs+M-Nb)$ for each $l\le s\le
n$, and then $\phi_l\le\dots\le\phi_n$. Thus $\phi_l= \min\{\phi_i:
1\le i\le n\}$.

We mention that if $A$ is symmetric, then the conditions (i)--(iii)
in Theorem~\ref{nonnega-up} hold if and only if for some $2\le t\le
l$, the following (i$^\prime$)--(iv$^\prime$) hold:

(i$^\prime$) $a_{ii}=M$ for $1\le i\le t-1$,

(ii$^\prime$) all off-diagonal elements of $A$ in the first $t-1$
rows and columns are equal to $N$,

(iii$^\prime$) $\frac{r_k(A)}{r_i(A)}=b$ for $1\le i\le n$, $1\le
k\le t-1$ and $k\ne i$,

(iv$^\prime$) $m_t=\dots=m_n$.

\noindent If $t=2$, then (i$^\prime$) and (ii$^\prime$) are
equivalent to $a_{11}=M$ and the off-diagonal elements of $A$ in the
first row and column are equal to $N$, and (iii$^\prime$) is
equivalent to $r_2(A)=\dots=r_n(A)$, implying that the conditions
(i$^\prime$)--(iv$^\prime$) above are equivalent to
(i$^{\prime\prime}$)--(iii$^{\prime\prime}$):

(i$^{\prime\prime}$) $a_{11}=M$ and the off-diagonal elements of $A$
in the first row and column are equal to $N$,

(ii$^{\prime\prime}$) $r_2(A)=\dots=r_n(A)$,

(iii$^{\prime\prime}$) $m_2=\dots=m_n$.

\noindent Suppose that $3\le t\le l$. From (i$^\prime$) and
(ii$^\prime$), $r_1(A)=\dots=r_{t-1}(A)=M+N(n-1)$. Thus from (iii'),
$b=\frac{r_2(A)}{r_1(A)}=1$, implying that $r_1(A)=\dots=r_n(A)$,
which also satisfies (iv'). It follows that
(i$^\prime$)--(iv$^\prime$) is equivalent to $A=MI_n+(N-M)J_n$,
which also satisfies conditions
(i$^{\prime\prime}$)--(iii$^{\prime\prime}$).

\begin{Theorem}\label{nonnega-low}
Let $A=(a_{ij})$ be an $n\times n$ nonnegative matrix with average
$2$-row sums $m_1,\dots,m_n$, where $m_1\ge\dots\ge m_n$. Let $S$ be
the smallest diagonal element, and $T$ the smallest off-diagonal
element of $A$. Let $c=\min\left\{\frac{r_j(A)}{r_i(A)}:1\le i,j\le
n\right\}$. Let
\[
\psi_n=\frac{m_{n}+S-Tc+
\sqrt{(m_{n}-S+Tc)^2+4Tc\sum_{i=1}^{n-1}(m_i-m_n)}}{2}.
\]
Then $\rho(A)\ge\psi_n$. Moreover, if $A$ is irreducible, then
$\rho(A)=\psi_n$ if and only if $m_1=\dots=m_n$ or $T>0$ and for
some $2\le t\le n$, $A$ satisfies the following conditions:

(i) $a_{ii}=S$ for $1\le i\le t-1$,

(ii) $m_t=\dots=m_n$,

(iii) $a_{ik}=T$ and $\frac{r_k(A)}{r_i(A)}=c$ for $1\le i\le n$ ,
$1\le k\le t-1$ and $k\ne i$.
\end{Theorem}

\begin{Proof}
For convenience, let $r_i=r_i(A)$ for $1\le i\le n$. Note that
$m_n=a_{nn}+\sum_{j=1}^{n-1}a_{nj}\frac{r_j}{r_n}\ge a_{ii}\ge S$.
If $T=0$, then $\psi_n=m_n$, and thus the result follows immediately
from Lemma~\ref{m1}. Suppose in the following that $T>0$.

Let $U=\text{diag}(r_1x_1,\dots,r_{n-1}x_{n-1},1)$, where $x_i\ge 1$
is a variable to be determined later for $1\le i\le n-1$. Let
$B=U^{-1}AU$. Obviously,  $A$ and $B$ are unitary similar, and thus
have the same eigenvalues. Recall that $m_i=\sum_{k=1}^n
a_{ik}\frac{r_k}{r_i}$ for $1\le i\le n$. For $1\le i\le n-1$, since
$a_{ii}\ge S$, and $a_{ik}\ge T$, $\frac{r_k}{r_i}\ge c$ for $k\ne
i$, we have
\begin{eqnarray*}
   r_i(B)
&=&
   \frac{1}{x_i}\left(\sum_{1\le k\le n-1\atop k\ne
   i}a_{ik}\frac{r_k}{r_i}(x_k-1)+a_{ii}(x_i-1)+m_i\right)\\
&\ge&
   \frac{1}{x_i}\left(Tc\sum_{1\le k\le n-1\atop k\ne
   i}(x_k-1)+S(x_i-1)+m_i\right)\\
&=&
   \frac{1}{x_i}\left(Tc\sum_{k=1}^{n-1}(x_k-1)+(S-Tc)(x_i-1)+m_i\right)
\end{eqnarray*}
with equality if and only if (a) and (b) hold: (a) $x_i=1$ or
$a_{ii}=S$, (b) $x_k=1$ or $a_{ik}=T$ and $\frac{r_k}{r_i}=c$, where
$1\le k\le n-1$ and $k\ne i$. Similarly,
\[
r_n(B)=\sum_{k=1}^{n-1}a_{nk}\frac{r_k}{r_n}(x_k-1)+m_n\ge
Tc\sum_{k=1}^{n-1}(x_k-1)+m_n
\]
with equality if and only if (c) holds: (c) $x_k=1$ or $a_{nk}=T$
and $\frac{r_k}{r_n}=c$, where $1\le k\le n-1$.

From the expression of $\psi_n$, we have
\[
Tc\sum_{k=1}^{n-1}(m_k-m_n)=(\psi_n-m_n)(\psi_n-S+Tc).
\]
Note that $Tc>0$. If $\sum_{k=1}^{n-1}(m_k-m_n)>0$, then
$\psi_n>\frac{m_n+S-Tc+|m_n-S+Tc|}{2}\ge\frac{m_n+S-Tc-(m_n-S+Tc)}{2}=S-Tc$,
and if $m_1=\dots=m_n$, then since $m_n\ge S$, we have
$\psi_n=m_n>S-Tc$. It follows that $\psi_n-S+Tc>0$. For $1\le i\le
n-1$, let $x_i=1+\frac{m_i-m_n}{\psi_n-S+Tc}$. Thus $x_i\ge 1$, and
\[
Tc\sum_{k=1}^{n-1}(x_k-1)
=\frac{Tc\sum_{k=1}^{n-1}(m_k-m_n)}{\psi_n-S+Tc}=\psi_n-m_n.
\]
Thus for $1\le i\le n-1$,
\begin{eqnarray*}
   r_i(B)
&\ge&
   \frac{1}{x_i}\left(Tc\sum_{k=1}^{n-1}(x_k-1)+(S-Tc)(x_i-1)+m_i\right)\\
&=&
   \frac{(\psi_n-m_n)+(S-Tc)\cdot\frac{m_i-m_n}{\psi_n-S+Tc}+m_i}{1+\frac{m_i-m_n}{\psi_n-S+Tc}}\\
&=&
   \psi_n,
\end{eqnarray*}
and
\[
r_n(B)\ge Tc\sum_{k=1}^{n-1}(x_k-1)+m_n=(\psi_n-m_n)+m_n=\psi_n.
\]
Hence by Lemma~\ref{r-min-max}, $\rho(A)=\rho(B)\ge\min_{1\le i\le
n}r_i(B)\ge\psi_n$.

Now suppose that $A$ is irreducible. Then $B$ is also irreducible.

Suppose that $\rho(A)=\psi_n$. Then $\rho(B)=\min_{1\le i\le
n}r_i(B)=\psi_n$, which, by Lemma~\ref{r-min-max}, implies that
$r_1(B)=\dots=r_n(B)=\psi_n$, and thus from the above arguments, (a)
and (b) for $1\le i\le n-1$ and (c) hold. If $m_1=m_n$, then
obviously $m_1=\dots=m_n$. Assume that $m_1>m_n$. Let $t$ be the
smallest integer such that $m_t=m_n$, where $2\le t\le n$. Then for
$1\le i\le t-1$, $m_i>m_n$, implying that $x_i>1$. Now (i)--(iii)
follow from (a), (b) for $1\le i\le n-1$ and (c).

Conversely, if $m_1=\dots=m_n$, then by Lemma~\ref{m1},
$\rho(A)=m_n=\psi_n$. If (i)--(iii) hold, then (a), (b) for $1\le
i\le n-1$ and (c) hold, implying that $r_i(B)=\psi_n$ for $1\le i\le
n$, and thus by Lemma~\ref{r-min-max}, $\rho(A)=\rho(B)=\psi_n$.
\end{Proof}

We mention that if $A$ is symmetric, then the conditions (i)--(iii)
in Theorem~\ref{nonnega-low} hold if and only if
(i$^{\prime\prime}$)--(iii$^{\prime\prime}$) hold:

(i$^{\prime\prime}$) $a_{11}=S$ and the off-diagonal elements of $A$
in the first row and column are equal to $T$,

(ii$^{\prime\prime}$) $r_2(A)=\dots=r_n(A)$,

(iii$^{\prime\prime}$) $m_2=\dots=m_n$.

To compare the above results with those in~\cite{DZ}, we listed the
corresponding results of~\cite{DZ} as follows.


\begin{Lemma}\label{A-up-Duan}\cite{DZ}
Let $A=(a_{ij})$ be an $n\times n$ nonnegative matrix with row sums
$r_1,\dots,r_n$, where  $r_1\ge\dots\ge r_n$. Let $M$ be the largest
diagonal element, and  $N$ the largest off-diagonal element of $A$.
Suppose that $N>0$. For $1\le l\le n$, let
\[
\Phi_l=\frac{r_{l}+M-N+\sqrt{(r_{l}-M+N)^2
+4N\sum_{i=1}^{l-1}(r_i-r_l)}}{2}.
\]
Then $\rho(A)\le\Phi_l$ for $1\le l\le n$. Moreover, if $A$ is
irreducible, then $\rho(A)=\Phi_l$ if and only if $r_1=\dots=r_n$ or
for some $2\le t\le l$, $A$ satisfies the following conditions:

(i) $a_{ii}=M$ for $1\le i\le t-1$,

(ii) $r_{t}=\dots=r_n$,

(iii) $a_{ik}=N$ for $1\le i\le n$, $1\le k\le t-1$ and $k\ne i$.
\end{Lemma}

\begin{Lemma}\label{A-low-Duan}\cite{DZ}
Let $A=(a_{ij})$ be an $n\times n$ nonnegative matrix with row sums
$r_1,\dots,r_n$, where $r_1\ge\dots\ge r_n$. Let $S$ be the smallest
diagonal element, and $T$ the smallest off-diagonal element of $A$.
Let
\[
\Psi_n= \frac{r_{n}+S-T+\sqrt{(r_{n}-S+T)^2
+4T\sum_{i=1}^{n-1}(r_i-r_n)}}{2}.
\]
Then $\rho(A)\ge\Psi_n$. Moreover, if $A$ is irreducible, then
$\rho(A)=\Psi_n$ if and only if $r_1=\dots=r_n$ or $T>0$, and for
some $2\le t\le n$, $A$ satisfies the following conditions:

(i) $a_{ii}=S$ for $1\le i\le t-1$,

(ii) $r_{t}=\dots=r_n$,

(iii) $a_{ik}=T$ for $1\le i\le n$ , $1\le k\le t-1$ and $k\ne i$.
\end{Lemma}

Consider
\[
A_1=\left(
\begin{array}{cccc}
0 & 1 & 1 & 1\\
1 & 0 & 2 & 2\\
1 & 2 & 0 & 2\\
1 & 2 & 2 & 0
\end{array}\right).
\]
In notations of Theorem~\ref{nonnega-up}, $m_1=5$,
$m_2=m_3=m_4=\frac{23}{5}$, $M=0$, $N=2$ and $b=\frac{5}{3}$,
implying that $\phi_1=5$, $\phi_2=4.7647$, $\phi_3=4.9230$ and
$\phi_4=5.0757$. It is easily seen that $A$ is unitary similar to
\[
A_1'=\left(
\begin{array}{cccc}
0 & 2 & 2 & 1\\
2 & 0 & 2 & 1\\
2 & 2 & 0 & 1\\
1 & 1 & 1 & 0
\end{array}\right),
\]
and thus $\rho(A_1)=\rho(A_1')$. In notations of
Lemma~\ref{A-up-Duan}, for $A'$, we have $r_1=r_2=r_3=5$, $r_4=3$,
$M=0$ and $N=2$, implying that $\Phi_1=\Phi_2=\Phi_3=5$ and
$\Phi_4=4.7720$. By direct check,
$\rho(A_1)=\rho(A_1')=4.6458<\min\{\phi_i:1\le i\le
4\}<\min\{\Phi_i:1\le i\le 4\}$.

Now consider
\[
A_2=\left(
\begin{array}{cccc}
0 & 0 & 0 & 1\\
0 & 0 & 0 & 1\\
0 & 0 & 0 & 1\\
1 & 1 & 1 & 0
\end{array}\right).
\]
In notations of Theorem~\ref{nonnega-up}, $m_1=m_2=m_3=3$, $m_4=1$,
$M=0$, $N=1$ and $b=3$, implying that $\phi_1=\phi_2=\phi_3=3$ and
$\phi_4=3.6904$. It is easily seen that $B$ is unitary similar to
\[
A_2'=\left(
\begin{array}{cccc}
0 & 1 & 1 & 1\\
1 & 0 & 0 & 0\\
1 & 0 & 0 & 0\\
1 & 0 & 0 & 0
\end{array}\right),
\]
and thus $\rho(A_2)=\rho(A_2')$. In notations of
Lemma~\ref{A-up-Duan}, for $B'$, we have $r_1=3$, $r_2=r_3=r_4=1$,
$M=0$ and $N=1$, implying that $\Phi_1=3$, $\Phi_2=1.732$,
$\Phi_3=2.236$ and $\Phi_4=2.6458$. By direct check,
$\rho(A_2)=\rho(A_2')=1.732=\min\{\Phi_i:1\le i\le
4\}<\min\{\phi_i:1\le i\le 4\}$.

The above  examples show that the upper  bound in Theorem~\ref{nonnega-up} and that in Lemma~\ref{A-up-Duan} are incomparable.

For $A_1$, in notations of Theorem~\ref{nonnega-low}, $m_1=5$,
$m_2=m_3=m_4=\frac{23}{5}$, $S=0$, $T=1$ and $c=\frac{3}{5}$,
implying that $\psi_4=4.6458$. For $A_1'$, in notations of
Lemma~\ref{A-low-Duan}, $r_1=r_2=r_3=5$, $r_4=3$, $S=0$ and $T=1$,
implying that $\Psi_4=4.1623$. Thus
$\rho(A_1)=\rho(A_1')=4.6458=\psi_4>\Psi_4$.

Consider
\[
A_3=\left(
\begin{array}{cccc}
0 & 2 & 2 & 4\\
2 & 0 & 2 & 2\\
2 & 2 & 0 & 2\\
1.9 & 1.9 & 1.9 & 0
\end{array}\right).
\]
In notations of Lemma~\ref{A-low-Duan}, we have $r_1=8$, $r_2=r_3=6$, $r_4=5.7$,
$S=0$ and $T=1.9$, implying that $\Psi_4=6.3665$. It is easily seen that $A_3$ is unitary similar to
\[
A_3'=\left(
\begin{array}{cccc}
0 & 1.9 & 1.9 & 1.9\\
2 & 0 & 2 & 2\\
2 & 2 & 0 & 2\\
4 & 2 & 2 & 0
\end{array}\right),
\]
and thus $\rho(A_3)=\rho(A_3')$. In notations of Theorem~\ref{nonnega-low}, for $A_3'$, we have $m_1=\frac{20}{3}$,
$m_2=m_3=\frac{197}{30}$, $m_4=5.85$, $S=0$, $T=1.9$ and $c=0.7125$,
implying that $\psi_4=6.2506$. By direct check,
$\rho(A_3)=\rho(A_3')>\Psi_4>\psi_4$.

The above  examples show that the lower  bound in Theorem~\ref{nonnega-low} and that in Lemma~\ref{A-low-Duan} are incomparable.

\section{Spectral radius of adjacency matrix}

Let $G$ be an $n$-vertex graph without isolated vertices. Let
$V(G)=\{v_1,\dots,v_n\}$. For $1\le i\le n$, recall that
$m_i(A(G))=\frac{\sum_{v_iv_j\in E(G)}d_j}{d_i}$ is  the average
$2$-degree of vertex $v_i$ in $G$. Let $d_{\max}$ and $d_{\min}$ be
respectively the maximal and minimal degrees of $G$. The following
result for a connected graph $G$ has been given by Huang and
Weng~\cite{HW}.

\begin{Theorem} \label{A-up}
Let $G$ be a graph on $n\ge 2$ vertices without isolated vertices.
Let $m_1\ge\dots\ge m_n$ be the average $2$-degrees of $G$. Then for
$1\le l\le n$,
\[
\rho(A(G))\le\frac{m_l-\frac{d_{\max}}{d_{\min}}
+\sqrt{\left(m_l+\frac{d_{\max}}{d_{\min}}\right)^2
+4\frac{d_{\max}}{d_{\min}}\sum_{i=1}^{l-1}(m_i-m_l)}}{2}.
\]
Moreover, if $G$ is connected, then equality holds if and only if
either $m_1=\dots=m_n$ or $d_1=n-1>d_2=\dots=d_n$.
\end{Theorem}

\begin{Proof}
We apply Theorem~\ref{nonnega-up} to $A(G)$. Since $M=0$, $N=1$ and
$b=\frac{d_{\max}}{d_{\min}}$, we have the desired upper bound for
$\rho(A(G))$. Note that $A(G)$ is symmetric. If $G$ is connected,
then $A(G)$ is irreducible, and thus the upper bound is attained if
and only if either $m_1=\dots=m_n$ or $A(G)$ satisfies the following
conditions (a)--(c):

(a) the off-diagonal elements of $A(G)$ in the first row and column
are equal to $1$,

(b) $r_2(A(G))=\dots=r_n(A(G))$,

(c) $m_2=\dots=m_n$,

\noindent or equivalently, either $m_1=\dots=m_n$ or (if $m_1>m_n$,
then) $d_1=n-1>d_2=\dots=d_n$.
\end{Proof}

From previous Theorem~\ref{A-up}, we have the following consequence:
Let $G$ be a graph on $n\ge 2$ vertices without isolated vertices.
Let $m_1\ge\dots\ge m_n$ be the average $2$-degrees of $G$. Then for
$1\le l\le n$,
\[
\rho(A(G))\le\frac{m_l-\frac{d_{\max}}{d_{\min}}
+\sqrt{\left(m_l+\frac{d_{\max}}{d_{\min}}\right)^2
+4\frac{d_{\max}}{d_{\min}}(l-1)(m_1-m_l)}}{2}.
\]
Moreover, if $G$ is connected, then equality holds if and only if
$m_1=\dots=m_n$.

\section{Spectral radius of signless Laplacian matrix}

Let $G$ be an $n$-vertex graph with $V(G)=\{v_1,\dots,v_n\}$. For
$1\le i\le n$, it is easily seen that
\[
m_i(Q(G))=d_i+\frac{\sum_{v_iv_j\in E(G)}d_j}{d_i},
\]
which is called the signless Laplacian average $2$-degree of vertex
$v_i$ in $G$. Recall that $d_{\max}$ and $d_{\min}$ are respectively
the maximal and minimal degrees of $G$ defined in Section~$3$.

\begin{Theorem} \label{q-max}
Let $G$ be a graph on $n\ge 2$ vertices without isolated vertices.
Let $m_1\ge\dots\ge m_n$ be the signless Laplacian average
$2$-degrees of $G$. Then for $1\le l\le n$,
\[
\rho(Q(G))\le\frac{m_l+d_{\max}-\frac{d_{\max}}{d_{\min}}
+\sqrt{\left(m_l-d_{\max}+\frac{d_{\max}}{d_{\min}}\right)^2
+4\frac{d_{\max}}{d_{\min}}\sum_{i=1}^{l-1}(m_i-m_l)}}{2}.
\]
Moreover, if $G$ is connected, then equality holds if and only if
$m_1=\dots=m_n$ or $d_1=n-1>d_2=\dots=d_n$.
\end{Theorem}

\begin{Proof}
We apply Theorem~\ref{nonnega-up} to $Q(G)$. Since $M=d_{\max}$,
$N=1$ and $b=\frac{d_{\max}}{d_{\min}}$, we have the desired upper
bound for $\rho(Q(G))$. Note that $Q(G)$ is symmetric. If $G$ is
connected, then $Q(G)$ is irreducible, and thus the upper bound is
attained if and only if either $m_1=\dots=m_n$ or $Q(G)=(q_{ij})$
satisfies the following conditions (a)--(c):

(a) $q_{11}=d_{\max}$, and the off-diagonal elements of $Q(G)$ in
the first row and column are equal to $1$,

(b) $r_2(Q(G))=\dots=r_n(Q(G))$,

(c) $m_2=\dots=m_n$,

\noindent or equivalently, either $m_1=\dots=m_n$ or (if $m_1>m_n$,
then) $d_1=n-1>d_2=\dots=d_n$.
\end{Proof}

\section{Spectral radius of  distance matrix}

Let $G$ be an $n$-vertex connected graph with
$V(G)=\{v_1,\dots,v_n\}$. For $1\le i\le n$, it is easily seen that
$m_i(D(G))=\frac{\sum_{j=1}^n d_{ij}D_j}{D_i}$, which is called the
average $2$-transmission of vertex $v_i$ in $G$. Let $D_{\max}$ and
$D_{\min}$ be respectively the maximal and minimal transmissions of
$G$.

\begin{Theorem} \label{D-up}
Let $G$ be a connected graph on $n\ge 2$ vertices. Let
$m_1\ge\dots\ge m_n$ be the average $2$-transmissions of $G$. For
$1\le l\le n$,
\[
\rho(D(G))\le \frac{m_{l}-\mathcal{D}\frac{D_{\max}}{D_{\min}}+
\sqrt{\left(m_{l}+\mathcal{D}\frac{D_{\max}}{D_{\min}}\right)^2
+4\mathcal{D}\frac{D_{\max}}{D_{\min}}\sum_{i=1}^{l-1}(m_i-m_l)}}{2}.
\]
with equality if and only if $m_1=\dots=m_n$.
\end{Theorem}

\begin{Proof}
We apply Theorem~\ref{nonnega-up} to $D(G)$. Since $M=0$,
$N=\mathcal{D}$ (the diameter of $G$) and
$b=\frac{D_{\max}}{D_{\min}}$, the upper bound follows with equality
if and only if either $m_1=\dots=m_n$ or (if $m_1>m_n$, then) $D(G)$
satisfies the following conditions (a)--(c):

(a) the off-diagonal elements of $D(G)$ in the first row and column
are equal to $\mathcal{D}$,

(b) $r_2(D(G))=\dots=r_n(D(G))$,

(c) $m_2=\dots=m_n$.

\noindent Since there is at least one $1$ in every row of $D(G)$,
(a) implies that $\mathcal{D}=1$, and thus $G$ is the $n$-vertex
complete graph, a contradiction for the latter case. Thus the upper
bound is attained if and only if $m_1=\dots=m_n$.
\end{Proof}

\begin{Theorem} \label{D-low}
Let $G$ be a connected graph on $n\ge 2$ vertices. Let
$m_1\ge\dots\ge m_n$ be the average $2$-transmissions of $G$. Then
\[
\rho(D(G))\ge\frac{m_n-\frac{D_{\min}}{D_{\max}}+
\sqrt{\left(m_n+\frac{D_{\min}}{D_{\max}}\right)^2
+4\frac{D_{\min}}{D_{\max}}\sum_{i=1}^{n-1}(m_i-m_n)}}{2}
\]
with equality if and only if $m_1=\dots=m_n$ or
$D_1=n-1<D_2=\dots=D_n$.
\end{Theorem}

\begin{Proof}
We apply Theorem~\ref{nonnega-low} to $D(G)$. Since $S=0$, $T=1$ and
$c=\frac{D_{\min}}{D_{\max}}$, we have the desired lower bound for
$\rho(D(G))$, which is attained if and only if either
$m_1=\dots=m_n$ or (if $m_1>m_n$, then) $D(G)$ satisfies the
following conditions (a)--(c):

(a) the off-diagonal elements of $D(G)$ in the first row and column
are equal to $1$,

(b) $r_2(D(G))=\dots=r_n(D(G))$,

(c) $m_2=\dots=m_n$,

\noindent or equivalently, either $m_1=\dots=m_n$ or (if $m_1>m_n$,
then) $D_1=n-1<D_2=\dots=D_n$.
\end{Proof}

\section{Spectral radius of  distance signless Laplacian matrix}

Let $G$ be an $n$-vertex connected graph with
$V(G)=\{v_1,\dots,v_n\}$. For $1\le i\le n$, it is easily seen that
\[
m_i(DQ(G))=D_i+\frac{\sum_{1\le j\le n\atop j\ne i}d_{ij}D_j}{D_i},
\]
which is called the signless Laplacian average $2$-transmission of
vertex $v_i$ in $G$. Recall that $D_{\max}$ and $D_{\min}$ are
respectively the maximal and minimal transmissions of $G$ defined in
Section~$5$.

Let $m_1\ge\dots\ge m_n$ be the signless Laplacian average
$2$-transmissions of $G$. We apply Theorem~\ref{nonnega-up} to
$DQ(G)$. Since $M=D_{\max}$, $N=\mathcal{D}$ (the diameter of $G$)
and $b=\frac{D_{\max}}{D_{\min}}$. By similar discussion as for the
distance matrix in Section~$5$, for $1\le l\le n$,
\[
\rho(DQ(G))\le
\frac{m_{l}+D_{\max}-\mathcal{D}\frac{D_{\max}}{D_{\min}}+
\sqrt{\left(m_{l}-D_{\max}+\mathcal{D}\frac{D_{\max}}{D_{\min}}\right)^2
+4\mathcal{D}\frac{D_{\max}}{D_{\min}}\sum_{i=1}^{l-1}(m_i-m_l)}}{2}
\]
with equality if and only if $m_1=\dots=m_n$.

\begin{Theorem}\label{DQ-low}
Let $G$ be a connected graph on $n\ge 2$ vertices. Let
$m_1\ge\dots\ge m_n$ be the signless Laplacian average
$2$-transmissions of $G$. Then
\[
\rho(DQ(G))\ge\frac{m_n+D_{\min}-\frac{D_{\min}}{D_{\max}}+
\sqrt{\left(m_n-D_{\min}+\frac{D_{\min}}{D_{\max}}\right)^2
+4\frac{D_{\min}}{D_{\max}}\sum_{i=1}^{n-1}(m_i-m_n)}}{2}
\]
with equality if and only if $m_1=\dots=m_n$ or
$D_1=n-1<D_2=\dots=D_n$.
\end{Theorem}

\begin{Proof}
We apply Theorem~\ref{nonnega-low} to $DQ(G)$. Since $S=D_{\min}$,
$T=1$ and $c=\frac{D_{\min}}{D_{\max}}$, we have the desired lower
bound for $\rho(DQ(G))$, and by similar arguments as in the proof of
Theorem~\ref{D-low}, it is attained if and only if either
$m_1=\dots=m_n$ or $D_1=n-1<D_2=\dots=D_n$.
\end{Proof}

\section{Spectral radius of reciprocal distance matrix}

Let $G$ be an $n$-vertex connected graph with
$V(G)=\{v_1,\dots,v_n\}$. For $1\le i\le n$, let
$R_i=r_i(R(G))=\sum_{1\le j\le n\atop j\ne i}\frac{1}{d_{ij}}$, and
then it is easily seen that
\[
m_i(R(G))=\frac{\sum_{1\le j\le n\atop j\ne
i}\frac{1}{d_{ij}}R_j}{R_i},
\]
which is called the  average $2$-reciprocal transmission of vertex
$v_i$ in $G$. Let $R_{\max}=\max\{R_i:1\le i\le n\}$ and
$R_{\min}=\min\{R_i:1\le i\le n\}$.

\begin{Theorem}\label{R-up}
Let $G$ be a connected graph on $n\ge 2$ vertices. Let
$m_1\ge\dots\ge m_n$ be the average $2$-reciprocal transmissions of
$G$. For $1\le l\le n$,
\[
\rho(R(G))\le\frac{m_l-\frac{R_{\max}}{R_{\min}}
+\sqrt{\left(m_l+\frac{R_{\max}}{R_{\min}}\right)^2
+4\frac{R_{\max}}{R_{\min}}\sum_{i=1}^{l-1}(m_i-m_l)}}{2}
\]
with equality if and only if $m_1=\dots=m_n$ or
$R_1=n-1>R_2=\dots=R_n$.
\end{Theorem}

\begin{Proof}
We apply Theorem~\ref{nonnega-up} to $R(G)$. Since $M=0$, $N=1$ and
$b=\frac{R_{\max}}{R_{\min}}$, we have the desired upper bound for
$\rho(R(G))$, which is attained if and only if either
$m_1=\dots=m_n$ or (if $m_1>m_n$, then) $R(G)$ satisfies the
following conditions (a)--(c):

(a) the off-diagonal elements of $R(G)$ in the first row and column
are equal to $1$,

(b) $r_2(R(G))=\dots=r_n(R(G))$,

(c) $m_2=\dots=m_n$,

\noindent or equivalently, $m_1=\dots=m_n$ or
$R_1=n-1>R_2=\dots=R_n$.
\end{Proof}

\vspace{3mm}

\noindent {\bf Acknowledgement.} This work was supported by the
National Natural Science Foundation of China (No.~11071089) and the
Specialized Research Fund for the Doctoral Program of Higher
Education of China (No.~20124407110002).

\end{document}